\documentclass[reqno,onecolumn,oneside]{paper}
\usepackage[english]{babel}
\usepackage{enumerate,amsmath,amsthm,amsfonts,pifont,slashed,amssymb,ifsym,mathrsfs,graphicx,graphics,yfonts,calligra,
latexsym,txfonts
}
\usepackage[all]{xy}
\usepackage[sort,compress]{cite}
\usepackage[usenames,dvipsnames]{color}
\newtheorem{theorem}{Theorem}[section]
\newtheorem{proposition}[theorem]{Proposition}
\newtheorem{lemma}[theorem]{Lemma}
\newtheorem{corollary}[theorem]{Corollary}

\theoremstyle{definition}
\newtheorem{definition}[theorem]{Definition}
\newtheorem{exm}[theorem]{Example}

\theoremstyle{remark}

\newcommand{\lcm}{\operatorname{lcm}}
\newcommand{\bydef}{\mathrel{\mathop:}=}

\newcommand{\obj}{\operatorname{Obj}}

\newcommand{\Rad}{\operatorname{Rad}}

\newcommand{\Aut}{\operatorname{Aut}}
\newcommand{\MVf}{\mathcal{MV}_f}

\newcommand{\Set}{\mathcal{S}\!\operatorname{et}}

\newcommand{\lang}{\mathcal{L}}

\newcommand{\Q}{\mathbb{Q}}

\newcommand{\cat}{\mathcal}

\newcommand{\B}{\operatorname{B}}
\newcommand{\sk}{\operatorname{Sk}}
\newcommand{\clop}{\operatorname{Clop}}
\newcommand{\Max}{\operatorname{Max}}

\newcommand{\LL}{\operatorname{\text{\L}}}

\newcommand{\epi}{\operatorname{epi}}

\newcommand{\R}{\mathbb{R}}
\newcommand{\N}{\mathbb{N}}

\newcommand{\TMV}{{}^{\operatorname{MV}}\!\mathcal{T}\!\!\operatorname{op}}
\renewcommand{\Top}{\mathcal{T}\!\!\operatorname{op}}
\renewcommand{\Set}{\mathcal{S}\!\operatorname{et}}

\newcommand{\Age}{\operatorname{Age}}
\newcommand{\Flim}{\operatorname{Flim}}

\renewcommand{\phi}{\varphi}

\newcommand{\cou}{{^\leftarrow}}

\newcommand{\restr}{\!\!\upharpoonright}
\newcommand{\la}{\langle}
\newcommand{\ra}{\rangle}

\newcommand{\lto}{\longrightarrow}

\newcommand{\ov}{\overline}

\newcommand{\0}{\mathbf{0}}
\renewcommand{\1}{\mathbf{1}}
\renewcommand{\2}{\mathbf{2}}

\def\amslatex\slash{{\protect\AmS-\protect\LaTeX}}

\begin{document} 

\title{Fraïssé limit and Ramsey Theorem: the case of MV-algebras and a categorical generalization}

\author{Ciro Russo} 

\institution{Departamento de Matemática \\ Instituto de Matemática e Estatística \\ Universidade Federal da Bahia -- Brazil}

\maketitle
\date{}

\begin{abstract}
In this paper we describe the Fraïssé limit of finite MV-algebras and then prove that finite MV-algebras verify the Ramsey property. Then we show that MV-algebras are just a special case of a more general situation. In fact, under minumum conditions for the application of the Kechris-Pestov-Todorcevic correspondence, Ramsey property holds for a certain category of finite structures if and only if it holds for a completion subcategory of it.
\end{abstract}

\section{Introduction}
\label{intro}

Structural Ramsey theory was introduced by Nešetřil and Rödl \cite{neserodl} as a more general approach to Ramsey theory \cite{ramsey} by means of some category-theoretical language.

Since then, the Ramsey and the Dual Ramsey properties have been proved to hold or not for many classes of mathematical structures, including a number of relevant algebraic ones. Structural Ramsey Theory is, on its turn, connected with Fraïssé Theory \cite{fra}, and is receiving quite a lot of interest since the discovery of its connection with topological dynamics \cite{kpt}.

MV-algebras are the equivalent algebraic semantics of \L ukasiewicz propositional logic and are intimately connected to Boolean algebras. In fact, they can be seen as a sort of non-idempotent generalization of Boolean algebras. Moreover, many properties of Boolean algebras either hold for MV-algebras or have a suitable MV-algebraic version. Such a connection has also a topological counterpart: Stone duality extends to a fuzzy topological duality, between a class of MV-algebras and a suitable fuzzy generalization of Stone spaces \cite{rusfuz}, that shall be used in the present paper.

In this short note we shall first describe the Fraïssé limit of the class of finite MV-algebras as a basis of rational-valued fuzzy clopens of the Cantor set equipped with the fuzzy topology induced by the usual one of $\R$, then we will prove that the finite Ramsey property holds for the class of MV-algebras, by using the Dual Ramsey property for finite Stone spaces (finite sets, in fact) \cite{grarot} and the fuzzy topological duality for MV-algebras \cite{rusfuz}. Finally, we will show that the way Ramsey property for Boolean algebras is inherited by MV-algebras is actually an instance of a more general fact. 

The paper is meant to be as much self-contained as possible, and is organized as follows. In Sections \ref{sec:mv} and \ref{sec:mvtop}, we recall the basic definitions about MV-algebras and MV-topologies, along with their respective properties that will be used in the main results. Section \ref{sec:frlim} is devoted to Fraïssé theory for finite MV-algebras: we first show that the class of finite MV-algebras is a Fraïssé class, and then we prove that the MV-algebra defined in (\ref{A}) is the Fraïssé limit of such a class (Theorem \ref{frlim}).

The Dual Ramsey property for finite-valued finite Stone MV-spaces (Theorem \ref{dualramseymv}) and, as a corollary, the Ramsey property for finite MV-algebras (Corollary \ref{ramseymv}) are proved in Section \ref{sec:ramsey}. Last, in Section \ref{ramseycatsec}, we prove that, if $\cat K$ and $\cat K'$ are two Fraïssé classes such that $\cat K$ is a completion subcategory of $\cat K'$, then the Ramsey property holds for $\cat K'$ if and only if it holds for $\cat K$ (Theorem \ref{ramseycat}). From such a result Corollary \ref{ramseymv} follows again by using also the Kechris-Pestov-Todor\v{c}evi\'{c} correspondence \cite{kpt}.

\section{MV-algebras}
\label{sec:mv}

MV-algebras were introduced in 1958 by C. C. Chang \cite{cha1} as an algebraic semantics -- the equivalent algebraic semantics, in fact -- of \L ukasiewicz propositional calculus. Since then, they have been quite extensively investigated not only as logic-related algebras, but also for their connection with other research areas and mathematical strutures. Indeed, besides being a pretty natural generalization of Boolean algebras, as we already stated, they are categorically equivalent to lattice-ordered Abelian groups with a strong order unit. The latter structures, on their turn, are the $K_0$ groups of some approximately finite-dimensional (AF, for short) C$^*$-algebras, and the $K_0$ group of any AF C$^*$-algebra embeds in one of them \cite{mun}.

Moreover, the connection between MV-algebras and Boolean algebras makes MV-algebras a natural algebraic framework for fuzzy set theory \cite{bel} and, consequently, for fuzzy topology \cite{rusfuz,dlprarc,dlprsoco}, and it is worth mentioning their connection with semirings and tropical mathematics \cite{dng,dnrcom}, with Riesz spaces \cite{dinletretr}, and their applications \cite{dnrins}.

In the present section we recall the main definitions and properties of MV-algebras. For further insights about them, we refer the reader to \cite{mvbook} and \cite{mubook}.

\begin{definition}\label{mvalg}
An \emph{MV-algebra} is an algebra $\la A, \oplus, \neg, 0 \ra$ of type $(2,1,0)$ that satisfies the following equations
\begin{enumerate}
\item[(\L 1)]$x \oplus (y \oplus z) = (x \oplus y) \oplus z$;
\item[(\L 2)]$x \oplus y = y \oplus x$;
\item[(\L 3)]$x \oplus 0 = x$;
\item[(\L 4)]$\neg(\neg x) = x$;
\item[(\L 5)]$x \oplus \neg 0 = \neg 0$;
\item[(\L 6)]$\neg(\neg x \oplus y) \oplus y = \neg (\neg y \oplus x) \oplus x$.
\end{enumerate}
\end{definition}

On every MV-algebra it is possible to define another constant $1$ and two further operations as follows:
\begin{enumerate}
\item[]$1 = \neg 0$,
\item[]$x \odot y = \neg (\neg x \oplus \neg y)$,
\item[]$x \ominus y = x \odot \neg y$.
\end{enumerate}
The following properties follow immediately from the definitions
\begin{enumerate}
\item[(\L 7)]$\neg 1 = 0$,
\item[(\L 8)]$x \oplus y = \neg (\neg x \odot \neg y)$,
\item[(\L 5')]$x \oplus 1 = 1$ \quad (reformulation of (\L 5)),
\item[(\L 6')]$(x \ominus y) \oplus y = (y \ominus x) \oplus x$ \quad (reformulation of (\L 6)),
\item[(\L 9)]$x \oplus \neg x = 1$.
\end{enumerate}

MV-algebras are naturally equipped with an order relation defined as follows
\begin{equation}
x \leq y \quad \textrm{ if and only if } \quad \neg x \oplus y = 1.
\end{equation}
Moreover it is easy to verify that $\neg x \oplus y = 1$ is indeed equivalent to each of the following conditions
\begin{itemize}
\item $x \odot \neg y = 0$;
\item $y = x \oplus (y \ominus x)$;
\item there exists an element $z \in A$ such that $x \oplus z = y$.
\end{itemize}

The order relation also determines a structure of bounded distributive lattice on $A$, with $0$ and $1$ respectively bottom and top element, and $\vee$ and $\wedge$ defined as follows
\begin{eqnarray}
&& x \vee y = (x \odot \neg y) \oplus y = (x \ominus y) \oplus y, \nonumber \\
&& x \wedge y = \neg (\neg x \vee \neg y) = x \odot (\neg x \oplus y). \nonumber 
\end{eqnarray}
It is worth noticing that $\oplus$, $\odot$, and $\wedge$ distribute over any existing join and, analogously, $\oplus$, $\odot$ and $\vee$ distribute over any existing meet, in any MV-algebra $A$. In other words, for any $Y \subseteq A$ for which there exists $\bigvee Y$, for any $Z \subseteq A$ for which there exists $\bigwedge Z$, and for any $x \in A$, 
$$x \bullet \bigvee Y  = \bigvee_{y \in Y} (a \bullet y), \text{ for } \bullet \in \{\oplus, \odot, \wedge\},$$
$$x \bullet \bigwedge Z = \bigwedge_{z \in Z} (x \bullet z), \text{ for } \bullet \in \{\oplus, \odot, \vee\}.$$

For any MV-algebra $A$, for all $a \in A$, and for all $n < \omega$, we shall use the  following abbreviations: $a^n \bydef \underbrace{a \odot \cdots \odot a}_{n \textrm{ times}}$ and $na \bydef \underbrace{a \oplus \cdots \oplus a}_{n \textrm{ times}}$. 

The best-known example of MV-algebra is the real unit interval $[0,1]$, with the sum $x \oplus y \bydef \min\{x + y, 1\}$ and the involution $\neg x \bydef 1 - x$. The product is, then, defined by $x \odot y \bydef \max\{x + y - 1, 0\}$, and the lattice structure resulting is the natural totally ordered one. It is well-known that the MV-algebra $[0,1]$ generates, both as a variety and as a quasi-variety, the whole variety of MV-algebras and that \L ukasiewicz propositional calculus is complete w.r.t. to such a structure. We also recall that, for any non-empty set $X$, $[0,1]^X$ with pointwise defined operations is an MV-algebra as well, and that Boolean algebras are precisely the MV-algebras satisfying the equation $x \oplus x = x$, hence they form a subvariety of the variety of MV-algebras.

For any  positive natural number $n$ the MV-algebra $\LL_n$ is the $n+1$-element subalgebra of $[0,1]$:
$$\LL_n = \{k/n \mid k = 0, \ldots, n\}.$$
Obviously, $\LL_1$ is the two-element Boolean algebra $\2$.

We recall that an \emph{ideal} of an MV-algebra $A$ is a non-empty downward closed subset $I$ such that $a \oplus b \in I$ for all $a, b \in I$. It is known that any non-trivial MV-algebra has maximal ideals; more precisely, any proper ideal of an MV-algebra is contained in a maximal one. The set of all maximal ideals of $A$ is denoted by $\Max A$, the ideal obtained as the intersection of all maximal ideals of $A$ is called the \emph{radical} of $A$: $\Rad A \bydef \bigcap \Max A$. \emph{Semisimple} algebras, as usual, are defined as the subdirect products of simple algebras. However, in the theory of MV-algebras, they can be characterized as (non-trivial) algebras whose radical is $\{0\}$, and such a characterization is most often used as a definition. It is worth noticing that Boolean algebras are all semisimple MV-algebras. 

Ideals and congruences of an MV-algebra $A$ are in bijective correspondence. Indeed, for any congruence $\theta$, $0/\theta$ is an ideal and, conversely, for any ideal $I$, the relation $\theta_I$ defined by ``$a \theta_I b$ iff $d(a,b) \bydef (a \odot \neg b) \oplus (b \odot \neg a) \in I$'' is a congruence on $A$ --- it is, in fact, the only one for which the class of zero is equal to $I$. Therefore, in MV-algebras, the congruence whose corresponding ideal is $I$ is often denoted by $I$ itself, and the congruence classes and the quotient algebra are denoted, respectively, by $a/I$ (for $a \in A$) and $A/I$.

Obviously \emph{simple} MV-algebras, i.e. those algebras whose congruence lattice is the two-element chain, have no non-trivial ideals. It may be worth recalling that ideals and filters are, in MV-algebras as well as in Boolean algebras, in one-one correspondence to each other. More precisely, a subset $I$ of an MV-algebra (or a Boolean algebra) $A$ is an ideal if and only if $\neg I = \{\neg a \mid a \in I\}$ is a filter. Such a correspondece is an order isomorphism between the sets of ideals and filters, ordered by set inclusion, and therefore $I$ is a maximal ideal if and only if $\neg I$ is an ultrafilter. In this paper, following the tradition of MV-algebras, we shall deal mainly with (maximal) ideals, but all the results and constructions presented can be suitably reformulated in terms of (ultra)filters.

Besides the fact that Boolean algebras are MV-algebras, it must be mentioned that, in the theory of MV-algebras, the so-called Boolean elements of each algebra play an important role.

An element $a$ of an MV-algebra $A$ is called \emph{idempotent} or \emph{Boolean} if $a \oplus a = a$. Boolean elements of MV-algebras enjoy the following properties.
\begin{itemize}
\item For any $a \in A$, $a \oplus a = a$ iff $a \odot a = a$.
\item An element $a$ is Boolean iff $\neg a$ is Boolean.
\item If $a$ and $b$ are idempotent, then $a \oplus b$ and $a \odot b$ are idempotent as well; moreover we have $a \oplus b = a \vee b$, $a \odot b = a \wedge b$, $a \vee \neg a = 1$, and $a \wedge \neg a = 0$.
\item The set $\B(A) = \{a \in A \mid a \oplus a = a\}$ is a Boolean algebra, usually called the \emph{Boolean center} of $A$. In fact, $\B(A)$ is the largest MV-subalgebra of $A$ which is Boolean.
\item For any $a \in A$ and $b \in \B(A)$, $a \oplus b = a \vee b$ and $a \odot b = a \wedge b$.
\end{itemize}

\begin{definition}\label{arch}
Let $A$ be an MV-algebra. An element $a \in A$ is \emph{Archimedean} if it satisfies the following equivalent conditions:
\begin{enumerate}[(a)]
\item there exists a positive integer $n$ such that $na \in \B(A)$;
\item there exists a positive integer $n$ such that $\neg a \vee na = 1$;
\item there exists a positive integer $n$ such that $na = (n + 1)a$.
\end{enumerate}
$A$ is called \emph{hyperarchimedean} if all of its elements are Archimedean.
\end{definition}

It is well-known that any hyperarchimedean MV-algebra is semisimple while the converse is not true. We recall a well-known representation theorem for semisimple MV-algebras.

\begin{theorem}\cite{bel}\label{belrepr}
For any set $X$, the MV-algebra $[0,1]^X$ and all of its subalgebras are semisimple. Moreover, up to isomorphisms, all the semisimple MV-algebras are of this type. More precisely, every semisimple MV-algebra can be embedded in the MV-algebra of fuzzy subsets $[0,1]^{\Max A}$ of the set of maximal ideals of $A$.
\end{theorem}

\section{Semisimple MV-algebras and fuzzy topologies}
\label{sec:mvtop}

In \cite{rusfuz} the author presented a pair of contravariant functors between a category of fuzzy topological spaces, called \emph{MV-topological spaces} or, simply, \emph{MV-topologies}, and the one of MV-algebras. Such functors, when restricted to suitable subcategories, form a duality which on its turn coincide with the classical Stone duality when restricted to Boolean algebras and Stone spaces (which are MV-topological spaces themselves).

Since such a duality will play an important role in our main results, we shall briefly recall pertinent definitions and results from \cite{rusfuz} in this section.

As we are going to see, the category of MV-topological spaces is a full subcategory of the one of fuzzy topological spaces in the sense of C. L. Chang \cite{chal}, and most of the basic definitions and results about them are simple adaptations of the corresponding ones from \cite{chal} or readily derivable from them.

\begin{definition}\label{mvtop}
Let $X$ be a set, $A$ the MV-algebra $[0,1]^X$, and $\tau \subseteq A$. We say that $\la X, \tau\ra$ is an \emph{MV-topological space} if 
\begin{enumerate}[(i)]
\item $\0, \1 \in \tau$,
\item for any family $\{o_i\}_{i \in I}$ of elements of $\tau$, $\bigvee_{i \in I} o_i \in \tau$,
\end{enumerate}
and, for all $o_1, o_2 \in \tau$,
\begin{enumerate}[(i)]
\setcounter{enumi}{2}
\item $o_1 \odot o_2 \in \tau$,
\item $o_1 \oplus o_2 \in \tau$,
\item $o_1 \wedge o_2 \in \tau$.
\end{enumerate}
$\tau$ is also called an \emph{MV-topology} on $X$ and the elements of $\tau$ are the \emph{open MV-subsets} of $X$.
\end{definition}

\begin{definition}\label{skeleton}
If $\la X, \tau\ra$ is an MV-topology, then $\la X, j(\tau)\ra$ --- where $j(\tau) \bydef \tau \cap \{0,1\}^X = \tau \cap \B([0,1]^X)$ --- is both an MV-topology and a topology in the classical sense. The topological space $\la X, j(\tau)\ra$ will be called the \emph{skeleton space} of $\la X, \tau\ra$.
\end{definition}

Note that $j(\tau)$ is the greatest topology contained in $\tau$, and $j$ is exactly the ``skeleton topology functor'' denoted by $\sk$ in \cite{rusfuz}.

The following two examples of MV-topologies are relevant to the content of next sections. With reference to the first one, we recall that, for any MV-algebra $A$ and $M \in \Max A$, the algebra $A/M$ is isomorphic to a (unique) subalgebra of $[0,1]$. The canonical embedding of $A/M$ into $[0,1]$ shall be denoted by $\iota_M$.
\begin{exm}\label{mvtopexm}
\begin{enumerate}
\item[(i)] For any MV-algebra $A$, consider the MV-algebra homomorphism $\iota: a \in A \mapsto (\iota_M(a/M))_{M \in \Max A} \in [0,1]^{\Max A}$ (which happens to be an embedding if and only if $A$ is semisimple). The image of $\iota$ is a basis (of clopens) for an MV-topology on $\Max A$, called the \emph{maximal MV-spectrum} of $A$ and denoted by $\tau_A$. The skeleton space of $\tau_A$ is the Stone topology on the maximal ideals of the Boolean center of $A$.
\item[(ii)] Let $d: X \lto [0,+\infty[$ be a distance function on $X$. For any fuzzy point $\alpha$ of $X$, with support $x$, and any positive real number $r$, we define the \emph{open ball} of center $\alpha$ and radius $r$ as the fuzzy set $\beta_r(\alpha)$ identified by the membership function $\beta_r(\alpha)(y) = \left\{\begin{array}{ll} \alpha(x) & \textrm{if } d(x,y) < r \\ 0 & \textrm{if } d(x,y) \geq r \end{array}\right.$. Analogously, the \emph{closed ball} $\beta_r[\alpha]$ of center $\alpha$ and radius $r$ has membership function $\beta_r[\alpha](y) = \left\{\begin{array}{ll} \alpha(x) & \textrm{if } d(x,y) \leq r \\ 0 & \textrm{if } d(x,y) > r \end{array}\right.$. It is immediate to verify that the fuzzy subsets of $X$ that are joins of a family of open balls is an MV-topology on $X$ that is said to be \emph{induced} by $d$ (see \cite{low}). The skeleton topology of such an MV-topology is the usual topology on the reals.
\end{enumerate}
\end{exm}

The MV-topology of Example \ref{mvtopexm}(i) defines a contravariant functor from the category of MV-algebras with their homomorphisms to the one of MV-topological spaces with fuzzy continuous functions. As in the case of Stone duality, by associating to each MV-topology the MV-algebra of its clopens, with obvious action on the morphisms, we obtain a contravariant functor $\clop$ in the opposite direction.

A semisimple MV-algebra is called \emph{limit cut complete} if it contains the join of a certain type of lattice cuts. The category of limit cut complete MV-algebra is a reflective subcategory of the one of all MV-algebras, and a completion category of the one of semisimple MV-algebra. The following is the content of Theorem 4.9 and Corollary 4.10 of \cite{rusfuz}. 

\begin{theorem}[Duality theorem]\label{dual}
$\clop$ and $\Max$ form a duality between the categories of limit cut complete MV-algebras and the one of fuzzy Hausdorff, fuzzy compact, zero-dimensional MV-topological spaces. The restriction of that pair of functors to Boolean algebras and Stone spaces (which are full subcategories of the two above) yields precisely the classical Stone Duality.
\end{theorem}

\section{The Fraïssé limit of finite MV-algebras}
\label{sec:frlim}

Fraïssé theory was introduced in 1954 \cite{fra} and boasts applications in various areas of mathematics, including functional analysis, topological dynamics, and Ramsey theory (see, for example, \cite{lupini,kpt,nesetril}). 

In the present section, we will first recall the main definitions of \cite{fra} and Fraïssé's theorem. After that, we will show that finite MV-algebras form a Fraïssé class, and we shall describe the Fraïssé limit of such a class. 

\begin{definition}\cite{fra}
Let $\lang$ be a first order signature and $A$ a countable $\lang$-structure which is locally finite (i.e., all of its finitely generated substructures are finite). \begin{itemize}
\item The \emph{age of $A$}, $\Age(A)$, is the class of all finite structures which can be embedded in $A$.
\item $A$ is said to be \emph{ultrahomogeneous} iff every isomorphism between two of its finite substructures can be extended to an automorphism of $A$.
\end{itemize}
\end{definition}
A locally finite, countably infinite, ultrahomogeneous structure shall be called a \emph{Fraïssé structure}.

We also recall the following 
\begin{definition}\label{class}
Let $\cat K$ a class of finitely generated $\lang$-structures. We say that $\cat K$
\begin{enumerate}
\item[(i)] has the \emph{hereditary property}, or that it is \emph{hereditary}, iff it contains all finitely generated substructures of all of its structures;
\item[(ii)] is \emph{essentially countable} iff it contains countably many structures up to isomorphisms;
\item[(iii)] has the \emph{joint embedding property} iff, for any two of its structures $A$ and $B$, there exists $C \in \cat K$ such that both $A$ and $B$ are embeddable in $C$;
\item[(iv)] has the \emph{amalgamation property} iff, given $A, B, C \in \cat K$ and two embeddings $f: A \to B$ and $g: A \to C$, there exists $D \in \cat K$ and two embeddings $h: B \to D$ and $k: C \to D$ such that $h \circ f = k \circ g$.
\end{enumerate}
\end{definition}

\begin{theorem}[Fraïssé]\label{frath}
Let $\lang$ be a signature and $\cat K$ a non-empty, essentially countable, hereditary class of finitely generated $\lang$-structures satisfying the joint embedding property and the amalgamation property. Then there exists a (unique, up to isomorphisms) countable $L$-structure $A$ such that $A$ is ultrahomogeneous and $\cat K = \Age(A)$.
\end{theorem}

\begin{definition}\label{flim}
A class satisfying (i--iv) of Definition \ref{class} is called a \emph{Fraïssé class}. The structure $A$ of Theorem \ref{frath} is a Fraïssé structure called the \emph{Fraïssé limit of $\cat K$} and denoted by $\Flim(\cat K)$.
\end{definition}

Now we are ready to show that finite MV-algebras form a Fraïssé class.

\begin{proposition}\cite[Proposition 3.6.5]{mvbook}\label{finitemv}
An MV-algebra $A$ is finite if and only if it is isomorphic to a finite product of finite chains,
$$A \cong \LL_{n_1} \times \cdots \times \LL_{n_t}, \text{ for some positive natural numbers } n_1, \ldots, n_t.$$
Moreover, this representation is unique up to the ordering of factors.
\end{proposition}

\begin{corollary}\label{finiteemb}
An MV-algebra $A$ is finite if and only there exist $n, t < \omega$ such that $A$ is embeddable in $\LL_n^t$.
\end{corollary}
\begin{proof}
The ``if'' part is trivial. The other one readily follows from Proposition \ref{finitemv} if one takes $n = \lcm(n_1, \ldots, n_t)$.
\end{proof}

\begin{proposition}\cite[Corollary 7.8 (iii)]{mubook}\label{coprod}
The class of finite MV-algebras is closed under finite free products.
\end{proposition}

\begin{corollary}\label{jep}
The class of finite MV-algebras enjoys the joint embedding property.
\end{corollary}

\begin{proposition}
The class of finite MV-algebras has the amalgamation property.
\end{proposition}
\begin{proof}
It is known (see, e.g., \cite[Section 2.6]{mubook}) that MV-algebras do enjoy the amalgamation property. On the other hand, the amalgamating object of a V-formation $B \leftarrow A \to C$ of finite MV-algebras is obviously a homomorphic image of the free product (i.e., the coproduct) of $B$ and $C$, and the latter is finite by Proposition \ref{coprod}. 
\end{proof}

By Proposition \ref{finitemv} and Corollary \ref{finiteemb} it follows immediately that the class of finite MV-algebras is essentially countable. The fact that it is also hereditary is obvious. Then it follows that the class $\MVf$ of finite MV-algebras is indeed a Fraïssé class which obviously contains the one of finite Boolean algebras. On the other hand, it is well-known that the Fraïssé limit of finite Boolean algebras is the countable atomless Boolean algebra $B_\infty$.

We will now describe the Fraïssé limit of finite MV-algebras and show that its automorphism group is isomorphic to the one of $B_\infty$, which happens to be its Boolean center.

\begin{proposition}\label{boocen}
$B(\Flim(\MVf)) \cong B_\infty$.
\end{proposition}
\begin{proof}
The assertion follows immediately from \cite[Proposition 2.3]{kpt} and the remarks following it.
\end{proof}

Let us now consider the Cantor set $(C, \tau)$ as a subspace of $\R$ with the usual topology and the MV-topology $\tau'$ on $C$ induced by the MV-topology on $\R$ defined in Example \ref{mvtopexm}(ii).

\begin{lemma}\label{Alem}
Let, for all $q \in [0,1] \cap \Q$, 
$$\ov q_U: x \in C \mapsto \left\{\begin{array}{ll} q & \text{if $x \in U$} \\ 0 & \text{if $x \notin U$}\end{array} \right..$$
Then the set
\begin{equation}\label{A}
A = \left\{\bigvee_{i=1}^n \ov{q_i}_{U_i} \mid n < \omega \& \forall i \leq n(q_i \in [0,1] \cap \Q \ \& \ U_i \in \clop\tau)\right\}
\end{equation}
is an MV-subalgebra of $[0,1]^C$ and a clopen basis for $\tau'$.
\end{lemma}
\begin{proof}
$A$ is obviously contained in $\tau'$ and both $\0 = \ov 0_C$ and $\1 = \ov 1 _C$ belong to $A$. For all $q \in [0,1] \cap \Q$ and $U \in \clop\tau$, $\neg \ov q_U = \ov 1_{U^c} \vee \ov{1-q}_U \in A$. Let $q, r \in [0,1] \cap \Q$ and $U, V \in \clop\tau$; then, for all $x \in C$,
$$(\ov q_U \wedge \ov r_V)(x) = \left\{\begin{array}{ll} 0 & \text{if $x \in (U \cap V)^c$} \\ q \wedge r & \text{if $x \in U \cap V$}\end{array}\right. = \ov{q \wedge r}_{U\cap V}.$$ 
Now, for $n > 1$, assume that $\neg \left(\bigvee_{i=1}^{n-1} \ov{q_i}_{U_i}\right) \in A$, i.e., there exist $V_1, \ldots, V_m \in \clop\tau$ and $r_1, \ldots, r_m \in [0,1]\cap \Q$ such that $\neg\left(\bigvee_{i=1}^{n-1} \ov{q_i}_{U_i}\right) = \bigvee_{j=1}^m \ov{r_j}_{V_j}$. We have
$$\begin{array}{l}
\neg\left(\bigvee_{i=1}^n \ov{q_i}_{U_i}\right) = \neg\left(\bigvee_{i=1}^{n-1} \ov{q_i}_{U_i}\right) \wedge \neg\ov{q_n}_{U_n} = \\
= \left(\bigvee_{j=1}^m \ov{r_j}_{V_j}\right) \wedge \neg\ov{q_n}_{U_n} = \bigvee_{j=1}^m (\ov{r_j}_{V_j} \wedge \neg\ov{q_n}_{U_n}) = \\
= \bigvee_{j=1}^m (\ov{r_j}_{V_j} \wedge (\ov 1_{U_n^c} \vee \ov{1-q_n}_{U_n})) = \bigvee_{j=1}^m ((\ov{r_j}_{V_j} \wedge \ov 1_{U_n^c}) \vee (\ov{r_j}_{V_j} \wedge \ov{1-q_n}_{U_n})) = \\
= \bigvee_{j=1}^m (\ov{r_j}_{V_j \cap U_n^c} \vee \ov{r_j \wedge (1-q_n)}_{V_j \cap U_n}) \in A,
\end{array}$$ 
and therefore it follows by induction on $n$ that $A$ is closed under $\neg$.

For what concerns the \L ukasiewicz sum, for $q, r \in [0,1] \cap \Q$, $U, V \in \clop\tau$, and $x \in C$, we have
$$(\ov q_U \oplus \ov r_V)(x) = \left\{\begin{array}{ll} 0 & \text{ if } x \in (U \cup V)^c \\
q \oplus r & \text{ if } x \in U \cap V \\
q & \text{ if } x \in U \setminus V \\
r & \text{ if } x \in V \setminus U
\end{array}\right.,$$
whence $\ov q_U \oplus \ov r_V = \ov{q \oplus r}_{U \cap V} \vee \ov q_{U \setminus V} \vee \ov r_{V \setminus U} \in A$.

It follows that $A \leq [0,1]^C$ and, by \cite[Proposition 3.3]{rusfuz}, $A$ is a base for an MV-topology on $C$.

Last, in order to prove that $A$ generates $\tau'$, it is sufficient to observe that, for any basis $X$ of the usual topology on $\R$, the set $\{\ov q_{U} \mid U \in X, q \in \Q\}$ is a basis for the MV-topology on $\R$ induced by the Euclidean distance.  
\end{proof}

\begin{lemma}\label{selfhomeo}
A map $f: C \to C$ is a self-homeomorphism w.r.t. $\tau'$ if and only if it is a self-homeomorphism w.r.t. $\tau$.
\end{lemma}
\begin{proof}
The right-to-left implication is obvious since $\tau \subseteq \tau'$ and the preimage of a crisp subset is crisp too.

Reciprocally, let us consider a basic clopen $\alpha = \bigvee_{i = 1}^n \ov{q_i}_{U_i}$; the support $U_i$ of each $\ov{q_i}_{U_i}$ is a clopen of $\tau$ and therefore so is its preimage under $f$. Therefore $\alpha \circ f = \bigvee_{i=1}^n \ov{q_i}_{f\cou(U_i)} \in A$, and an analogous argument holds for $f^{-1}$. The assertion follows.
\end{proof}

\begin{theorem}\label{frlim}
$\Flim(\MVf) \cong A$. 
\end{theorem}
\begin{proof}
It is obvious that every finite MV-algebra embeds in $A$ and that $A$ is countable because it is the union of countably many countable sets. Moreover, every finitely generated subalgebra of $A$ is trivially finite. So we only need to prove that $A$ is ultrahomogeneous.

By Proposition \ref{boocen} and the definition of $A$, $\Flim(\MVf)$ and $A$ have the same Boolean center. If $D$ and $E$ are finitely generated subalgebras of $A$ and $f: D \to E$ is an isomorphism, the restriction $f\restr: B(D) \to B(E)$ is an isomorphism between their respective Boolean centers; since $B_\infty$ is the Fraïssé limit of finite Boolean algebras, there exists an automorphism $g$ of $B_\infty \cong \clop\tau$ extending $f\restr$. By applying the Stone functor we obtain a self-homeomorphism $h = \Max\restr g$ of $(C, \tau)$ whose underlying map, by Lemma \ref{selfhomeo}, is a self-homeomorphism of the MV-topological space $(C,\tau')$. Then $\clop h$ is an automorphism of $A$ which extends $f$. 

It follows that $A$ is a countable and ultrahomogeneous MV-algebra containing isomorphic copies of each finite MV-algebra, whence it is, up to isomorphisms, the Fraïssé limit of the class of finite MV-algebras.
\end{proof}

\begin{corollary}\label{autiso}
$\Aut(\Flim(\MVf)) \cong \Aut(B_\infty)$.
\end{corollary}
\begin{proof}
It follows immediately from Lemma \ref{selfhomeo} and Theorem \ref{frlim}.
\end{proof}

\section{Ramsey Theorem for finite MV-algebras}
\label{sec:ramsey}

In the present section we shall prove the Ramsey property for the class of finite MV-algebras. Let us recall what it means for a category to have the (structural) Ramsey Property (see \cite{leeb,grlero}).
\begin{definition}
Let $\cat C$ be a category. $\cat C$ has the \emph{Ramsey property} iff for every natural number $r \geq 2$, and all objects $A, B \in \obj\cat C$ (with non-empty $\hom_{\cat C}(A,B)$), there exists $D \in \obj\cat C$, such that for every $r$-colouring $\Gamma: \hom_{\cat C}(A,D) \to r$ there exists a $\cat C$-morphism $f: B \to D$ such that the set
$\{f \circ g \mid g \in \hom_{\cat C}(A,B)\}$ is $\Gamma$-monochromatic.

$\cat C$ is said to satisfy the \emph{Dual Ramsey property} iff the opposite category $\cat C^{\operatorname{op}}$ has the Ramsey property.
\end{definition}

In the next three results we will show that the Dual Ramsey property for finite sets implies the Ramsey property for finite MV-algebras.

For any set $X$ of cardinality $m < \omega$ and for any $n < \omega$, we shall denote by $\tau_{X,n}$ the set $\LL_n^X$. Note that $\tau_{X,n}$ is the finest $(n+1)$-valued MV-topology on $X$. 

\begin{lemma}\label{homset}
Let $n < \omega$, $\la X, \tau_X \ra$ and $\la Y, \tau_Y\ra$ $(n+1)$-valued MV-topological spaces, and let us denote by $j(X)$ and $j(Y)$ the skeleton topological spaces $\la X, j(\tau_X)\ra$ and $\la Y, j(\tau_Y)\ra$ respectively.

Then $\hom_{\TMV}(X,Y) \subseteq \hom_{\Top}(j(X),j(Y))$. Moreover, the equality holds if $\tau_X = \tau_{X,n}$.   
\end{lemma}
\begin{proof}
For all $f \in \hom_{\TMV}(X,Y)$, the preimage under $f$ of a Boolean element of $\tau_Y$ is a Boolean element of $\tau_X$. Therefore the inclusion follows immediately from \cite[Proposition 3.7 (vii)]{dlprarc}. If $\tau_X = \tau_{X,n}$, any function from $X$ to $Y$ is obviously MV-continuous, whence the other inclusion follows.
\end{proof}

Before proving next result, we recall that epimorphisms in the categories $\TMV$, $\Top$, and $\Set$ are exactly the morphisms whose underlying map is surjective. 

\begin{theorem}\label{dualramseymv}
Dual Ramsey Theorem holds for finite-valued finite Stone MV-\-topo\-log\-i\-cal spaces.
\end{theorem}
\begin{proof}
Let $\la X, \tau_X\ra$ and $\la Y, \tau_Y\ra$ be finite-valued finite Stone MV-topological spaces, with $\left| Y \right| = k > m = \left| X \right|$, and let us denote by $\epi_{\TMV}(Y,X)$ the set of epimorphisms (i.e., surjective MV-continous maps) from $Y$ to $X$. $j(X)$ and $j(Y)$ are then Stone spaces of cardinalities $m$ and $k$ respectively, and we have that $\epi_{\TMV}(Y,X) \subseteq \epi_{\Top}(j(Y),j(X)) = \epi_{\Set}(Y,X)$, by Lemma \ref{homset}. The Dual Ramsey Theorem for finite sets or, equivalently, for finite Stone spaces, \cite{grarot, nesrodproc, carsim} guarantees that there exists a finite set $Z$ of cardinality, say, $N$, and an epimorphism $p: Z \to Y$ such that, for any $r \in \omega$ and for any $r$-colouring $\Gamma$ of $\epi_{\Set}(Z, X)$, the set $\{f \circ p \mid f \in \epi_{\Set}(Y,X)\}$ is $\Gamma$-monochromatic. 

Let $n_X$ and $n_Y$ be the least common denominators of the values (written as irreducible fractions) taken by the elements of $\tau_X$ and $\tau_Y$ respectively, and let $n = \lcm(n_X,n_Y)$. Consider the $(n+1)$-valued MV-Stone space $\la Z, \tau_{Z,n}\ra$. By Lemma \ref{homset}, we have:
$$\begin{array}{l} 
\epi_{\TMV}(\la Z, \tau_{Z,n}\ra, \la X, \tau_X\ra) = \epi_{\Set}(Z,X), \text{ and} \\
\epi_{\TMV}(\la Z, \tau_{Z,n}\ra, \la Y, \tau_Y\ra) = \epi_{\Set}(Z,Y).
\end{array}$$
Such equalities imply:
\begin{itemize}
\item $p \in \epi_{\TMV}(\la Z, \tau_{Z,n}\ra, \la Y, \tau_Y\ra),$
\item $\{f \circ p \mid f \in \epi_{\TMV}(Y,X)\} \subseteq \{f \circ p \mid f \in \epi_{\Set}(Y,X)\}$, and
\item any $r$-colouring of $\epi_{\TMV}(\la Z, \tau_{Z,n}\ra, \la X, \tau_X\ra)$ is an $r$-colouring of $\epi_{\Set}(Z,X)$.
\end{itemize}
It follows that $\{f \circ p \mid f \in \epi_{\TMV}(Y,X)\}$ is $\Gamma$-monochromatic. 
\end{proof}

\begin{corollary}\label{ramseymv}
Ramsey Theorem holds for finite MV-algebras.
\end{corollary}
\begin{proof}
It follows immediately from Theorem \ref{dualramseymv} by observing that finite MV-algebras are dual to finite-valued finite Stone MV-spaces in the duality of Theorem \ref{dual}. 
\end{proof}

\section{Transferring the Ramsey property between Fraïssé classes}
\label{ramseycatsec}

The validity of Ramsey property for finite MV-algebras strongly relies on the one for finite Boolean algebras, as we saw in the previous sections. It turns out that this is just a special case of a definitely more general fact, namely, that Ramsey property is preserved under categorical completions, given obvious additional hypotheses meant to ensure that the given categories are Fraïssé classes, as we show in Theorem \ref{ramseycat}.

Let us first recall the definition of extremely amenable group and the Kechris-Pestov-Todor\v{c}evi\'{c} correspondence.
\begin{definition}
A topological group $G$ is called extremely amenable iff every continuous action $G \curvearrowright K$ on a compact space $K$ has a fixed point, i.e., there exists $p \in K$ such that $a \cdot p = p$ for all $a \in G$.
\end{definition}

\begin{theorem}[Theorem 4.7 of \cite{kpt}]\label{kptth}
Let $G$ be a closed subgroup of $\Aut(\Q, <)$. Then the following are equivalent:
\begin{enumerate}
\item[(a)]$G$ is extremely amenable;
\item[(b)] $G = \Aut(A)$, where $A$ is the Fraïssé limit of a Fraïssé order class with the Ramsey property.
\end{enumerate}
\end{theorem}
Observe that a \emph{Fraïssé order class} is a Fraïssé class with a distinguished total order in its language.

We can now generalize Theorem \ref{dualramseymv}.
\begin{theorem}\label{ramseycat}
Let $\cat D$ be a category and $\cat C$ be a completion for $\cat D$, i.e., there exist functors $F: \cat C \to \cat D$ and $G: \cat D \to \cat C$ such that:
\begin{enumerate}
\item[(i)] $F$ is a full embedding,
\item[(ii)] $G$ is left inverse and left adjoint to $F$,
\item[(iii)] $G$ is faithful.
\end{enumerate} 
Moreover, assume that $\cat K$ and $\cat K'$ are Fraïssé classes with a distinguished total order, made of objects of $\cat C$ and $\cat D$ respectively, such that $G[\cat K'] = \cat K$.

Then $\cat K$ satisfies the Ramsey property if and only if so does $\cat K'$.
\end{theorem}
\begin{proof}
Let $C = \Flim(\cat K)$ and $D = \Flim(\cat K')$. From a category-theoretic viewpoint and up to isomorphisms, $C$ and $D$ are the colimit objects, in $\cat C$ and $\cat D$ respectively, of their respective ages, with any given diagrams of embeddings $d: I \to \cat C$ and $d': J \to \cat D$, where $I, J \subseteq \N$ are the (unbounded, naturally ordered) sets of possible cardinalities for the finite structures of $\cat C$ and $\cat D$ respectively (see also \cite{Kubis2007FrassSC,masulovic1}, and \cite{Bartovs2021TheWR} for categorical approaches to Fraïssé classes and limits). Therefore $G(D) = C$ because completions preserves colimits.

On the other hand, the hom-set restrictions of $G$ are also bijective by (ii) and (iii), whence the groups $\Aut C$ and $\Aut D$ are isomorphic. By Theorem \ref{kptth}, a Fraïssé order class satisfies the Ramsey property if and only if the automorphism group of its Fraïssé limit is extremely amenable. Then the assertion follows from the isomorphism between $\Aut C$ and $\Aut D$.
\end{proof}

As a consequence of Theorem \ref{ramseycat}, we can immediately obtain once again Corollary \ref{ramseymv} by simply observing that the category of Boolean algebras with Boolean embeddings as morphisms is a completion of the one of MV-algebras with MV-algebra embeddings, where the completion functor is the one that associates to each MV-algebra its own Boolean center and acts as the restriction on morphisms, and its right adjoint is simply the inclusion. In this case, both the classes at hand can easily be made into Fraïssé order ones. Indeed, all finite Boolean and MV-algebras are subdirect products of totally ordered (simple) algebras, and therefore the lexicographic order of the product of such factors is total and induces a total ordering on its subalgebras.

\bibliographystyle{alphadin}

\bibliography{cirobibtex}

\begin{thebibliography}{DLPR20b}


\providecommand{\url}[1]{\texttt{#1}}
\expandafter\ifx\csname urlstyle\endcsname\relax
  \providecommand{\doi}[1]{doi: #1}\else
  \providecommand{\doi}{doi: \begingroup \urlstyle{rm}\Url}\fi

\bibitem[BBBK21]{Bartovs2021TheWR}
\textsc{Bartoš}, A. ; \textsc{Bice}, T. ; \textsc{Barbosa}, K.~D.  ;
  \textsc{Kubiś}, W.:
\newblock The weak Ramsey property and extreme amenability.
\newblock {In: }\emph{arXiv preprint arXiv:2110.01694 [math.LO]}  (2021)

\bibitem[Bel86]{bel}
\textsc{Belluce}, L.~P.:
\newblock Semisimple algebras of infinite valued logic and bold fuzzy set
  theory.
\newblock {In: }\emph{Can. J. Math} 38 (1986), Nr. 6, S. 1356--1379

\bibitem[CDM00]{mvbook}
\textsc{Cignoli}, R. L.~O. ; \textsc{D'Ottaviano}, I. M.~L.  ;
  \textsc{Mundici}, D.:
\newblock \emph{Trends in Logic}. Bd.~7: {\emph{Algebraic Foundations of
  Many-valued Reasoning}}.
\newblock Kluwer Academic Publishers, Dordrecht, 2000

\bibitem[Cha58]{cha1}
\textsc{Chang}, C.~C.:
\newblock Algebraic analysis of many valued logic.
\newblock {In: }\emph{Trans. Amer. Math. Soc} 88 (1958), S. 467--490

\bibitem[Cha68]{chal}
\textsc{Chang}, C.~L.:
\newblock Fuzzy topological spaces.
\newblock {In: }\emph{J. Math. Anal. App} 24 (1968), Nr. 1, S. 182--190

\bibitem[CS84]{carsim}
\textsc{Carlson}, T.~J. ; \textsc{Simpson}, S.~G.:
\newblock A dual form of Ramsey Theorem.
\newblock {In: }\emph{Advances in Mathematics} 53 (1984), Nr. 3, S. 265--290

\bibitem[DLPR20a]{dlprarc}
\textsc{De~La~Pava}, L.~V. ; \textsc{Russo}, C.:
\newblock Compactness in {MV}-Topologies: {T}ychonoff Theorem and
  {Stone-\v{C}ech} Compactification.
\newblock {In: }\emph{Archive for Mathematical Logic} 59 (2020), S. 57--79

\bibitem[DLPR20b]{dlprsoco}
\textsc{De~La~Pava}, L.~V. ; \textsc{Russo}, C.:
\newblock MV-{A}lgebras as Sheaves of {$\ell$}-groups on Fuzzy Topological
  Spaces.
\newblock {In: }\emph{Soft Computing} 24 (2020), S. 8793--8804

\bibitem[DNG05]{dng}
\textsc{Di~Nola}, A. ; \textsc{Gerla}, B.:
\newblock Algebras of {{\L}ukasiewicz}'s Logic and their Semiring Reducts.
\newblock {In: }\emph{Idempotent Mathematics and Mathematical Physics -
  Contemp. Math 377}.
\newblock Amer. Math. Soc, 2005, S. 131--144

\bibitem[DNL96]{dinletretr}
\textsc{Di~Nola}, A. ; \textsc{Lettieri}, A.:
\newblock Coproduct {MV}-Algebras, nonstandard reals, and {R}iesz Spaces.
\newblock {In: }\emph{Journal of Algebra} 185 (1996), S. 605--620

\bibitem[DNR07]{dnrins}
\textsc{Di~Nola}, A. ; \textsc{Russo}, C.:
\newblock {{\L}ukasiewicz} {T}ransform and its application to compression and
  reconstruction of digital images.
\newblock {In: }\emph{Information Sciences} 177 (2007), S. 1481--1498

\bibitem[DNR13]{dnrcom}
\textsc{Di~Nola}, A. ; \textsc{Russo}, C.:
\newblock Semiring and semimodule issues in {MV}-algebras.
\newblock {In: }\emph{Communications in Algebra} 41 (2013), Nr. 3, S.
  1017--1048

\bibitem[Fra54]{fra}
\textsc{Fraïssé}, R.:
\newblock Sur l'extension aux relations de quelques proprietés des ordres.
\newblock {In: }\emph{Ann. Sci. École Norm. Sup.} 71 (1954), S. 363--388

\bibitem[GLL72]{grlero}
\textsc{Graham}, R.~L. ; \textsc{Leeb}, K.  ; \textsc{L.}, Rothschild~B.:
\newblock Ramsey Theorem for a class of categories.
\newblock {In: }\emph{Advances in Mathematics} 8 (1972), Nr. 3, S. 417--433

\bibitem[GR71]{grarot}
\textsc{Graham}, R.~L. ; \textsc{Rothschild}, B.~L.:
\newblock Ramsey's Theorem for $n$-Parameter Sets.
\newblock {In: }\emph{Transactions of the American Mathematical Society} 159
  (1971), S. 257--292

\bibitem[KPT05]{kpt}
\textsc{Kechris}, A.~S. ; \textsc{Pestov}, V.  ; \textsc{Todor\v{c}evi\'{c}},
  S.:
\newblock Fraïssé limits, {R}amsey theory and topological dynamics of
  automorphism groups.
\newblock {In: }\emph{Geom. Funct. Anal.} 15 (2005), Nr. 1, S. 106--189

\bibitem[Kub07]{Kubis2007FrassSC}
\textsc{Kubi{\'{s}}}, W.:
\newblock Fra{\"{i}}ss{\'{e}} sequences: category-theoretic approach to
  universal homogeneous structures.
\newblock {In: }\emph{Ann. Pure Appl. Log.} 165 (2007), S. 1755--1811

\bibitem[Lee73]{leeb}
\textsc{Leeb}, K.:
\newblock Vorlesungen über {P}ascaltheorie.
\newblock {In: }\emph{Arbeitsberichte des Inst. für math. Maschinen und
  Datenverarbeitung -- Friedrich Alexander Universität Erlangen Nürnberg} 6
  (1973), Nr. 7

\bibitem[Low76]{low}
\textsc{Lowen}, R.:
\newblock Fuzzy topological spaces and fuzzy compactness.
\newblock {In: }\emph{J. Math. Anal. Appl} 56 (1976), S. 621--633

\bibitem[Lup18]{lupini}
\textsc{Lupini}, M.:
\newblock Fraïssé limits in functional analysis.
\newblock {In: }\emph{Advances in Mathematics} 338 (2018), S. 93--174

\bibitem[Ma{\v{s}}20]{masulovic1}
\textsc{Ma{\v{s}}ulovi{\'{c}}}, D.:
\newblock The Kechris–Pestov–Todor{\v{c}}evi{\'{c}} Correspondence from the
  Point of View of Category Theory.
\newblock {In: }\emph{Applied Categorical Structures} 29 (2020), S. 141 -- 169

\bibitem[Mun86]{mun}
\textsc{Mundici}, D.:
\newblock Interpretation of {$\mathnormal{AF C}^*$}-algebras in
  {{\L}ukasiewicz} sentential calculus.
\newblock {In: }\emph{J. Functional Analysis} 65 (1986), S. 15--63

\bibitem[Mun11]{mubook}
\textsc{Mundici}, D.:
\newblock \emph{Trends in Logic}. Bd.~35: {\emph{Advanced {{\L}ukasiewicz}
  Calculus and {MV}-algebras}}.
\newblock Springer, 2011

\bibitem[Ne{\v{s}}05]{nesetril}
\textsc{Ne{\v{s}}et{\v{r}}il}, J.:
\newblock Ramsey Classes and Homogeneous Structures.
\newblock {In: }\emph{Combinatorics, Probability and Computing} 14 (2005), Nr.
  1–2, S. 171–189

\bibitem[NR77]{neserodl}
\textsc{Ne{\v{s}}et{\v{r}}il}, J. ; \textsc{Rödl}, V.:
\newblock Partitions of finite relational and set systems.
\newblock {In: }\emph{Journal of Combinatorial Theory, Series A} 22 (1977), Nr.
  3, S. 289--312

\bibitem[NR80]{nesrodproc}
\textsc{Ne{\v{s}}et{\v{r}}il}, J. ; \textsc{Rödl}, V.:
\newblock {D}ual {R}amsey {T}heorem.
\newblock {In: }\textsc{Frolik}, Z. (Hrsg.): \emph{Proc. 8th Winter School on
  Abstract Analysis}.
\newblock Prague : Czechoslovak Academy of Sciences, 1980, S. 121--123

\bibitem[Ram30]{ramsey}
\textsc{Ramsey}, F.~P.:
\newblock On a Problem of Formal Logic.
\newblock {In: }\emph{Proceedings of the London Mathematical Society} s2-30
  (1930), Nr. 1, S. 264--286

\bibitem[Rus16]{rusfuz}
\textsc{Russo}, C.:
\newblock An extension of {S}tone {D}uality to fuzzy topologies and
  {MV}-algebras.
\newblock {In: }\emph{Fuzzy Sets and Systems}  (2016), Nr. 303, S. 80--96

\end{thebibliography}

\end{document}